\documentclass[aos,MSNbibl,seceqn,citesort,dvips]{arximspdf}
\usepackage{graphicx}

\doi{10.1214/11-AOS910}
\volume{39}
\issue{5}
\pubyear{2011}
\firstpage{2533}
\lastpage{2556}

\makeatletter

\def\bsuffix #1{#1}

\newtheorem{theorem}{Theorem}
\newtheorem{prp}{Proposition}
\newtheorem{cor}{Corollary}

\newcommand{\goto}{\rightarrow}

\newcommand{\vct}[1]{\mathbf{#1}}
\newcommand{\mtx}[1]{\mathbf{#1}}

\newcommand{\eps}{\varepsilon}

\newcommand{\cS}{\mathcal{S}}
\newcommand{\cN}{\mathcal{N}}
\newcommand{\cJ}{\mathcal{J}}

\newcommand{\bA}{\mathbf{A}}
\newcommand{\bX}{\mathbf{X}}
\newcommand{\bI}{\mathbf{I}}
\newcommand{\by}{\mathbf{y}}
\newcommand{\bz}{\mathbf{z}}
\newcommand{\bbeta}{{\bolds\beta}}

\newcommand{\bbN}{\mathbb{N}}
\newcommand{\bbR}{\mathbb{R}}

\makeatother

\begin{document}
\begin{frontmatter}

\title{Global testing under sparse alternatives:
ANOVA, multiple comparisons and the higher criticism\thanksref{T1}}
\runtitle{Global testing under sparse alternatives}

\thankstext{T1}{Supported in part by an ONR Grant N00014-09-1-0258.}

\begin{aug}
\author[A]{\fnms{Ery} \snm{Arias-Castro}\corref{}\ead[label=e1]{eariasca@ucsd.edu}},
\author[B]{\fnms{Emmanuel J.} \snm{Cand\`es}\ead[label=e2]{candes@stanford.edu}}
\and
\author[C]{\fnms{Yaniv} \snm{Plan}\ead[label=e3]{plan@caltech.edu}}
\runauthor{E. Arias-Castro, E. J. Cand\`es and Y. Plan}
\affiliation{University of California, San Diego, Stanford
University\break
and California Institute of Technology}
\address[A]{E. Arias-Castro\\
Department of Mathematics\\
University of California, San Diego\\
9500 Gilman Drive\\
San Diego, California 92093-0112\\
USA\\
\printead{e1}}
\address[B]{E. J. Cand\`es\\
Department of Statistics\\
Stanford University\\
390 Serra Mall\\
Stanford, California 94305-4065\\
USA\\
\printead{e2}}
\address[C]{Y. Plan\\
Department of Applied\\
\quad and Computational Mathematics\\
California Institute of Technology\\
300 Firestone, Mail Code 217-50\\
Pasadena, California 91125\\
USA\\
\printead{e3}} 
\end{aug}

\received{\smonth{7} \syear{2010}}
\revised{\smonth{4} \syear{2011}}

%
\begin{abstract}
Testing for the significance of a subset of regression coefficients
in a linear model, a staple of statistical analysis, goes back at
least to the work of Fisher who introduced the analysis of variance
(ANOVA).
We study this problem under the assumption that the coefficient
vector is sparse, a common situation in modern high-dimensional
settings. Suppose we have $p$ covariates and that under the
alternative, the response only depends upon the order of
$p^{1-\alpha}$ of those, $0 \le\alpha\le1$. Under moderate
sparsity levels, that is, $0 \le\alpha\le1/2$, we show that ANOVA is
essentially optimal under some conditions on the design. This is no
longer the case under strong sparsity constraints, that is, $\alpha>
1/2$. In such settings, a multiple comparison procedure is often
preferred and we establish its optimality when $\alpha\geq3/4$.
However, these two very popular methods are suboptimal, and
sometimes powerless, under moderately strong sparsity where $1/2 <
\alpha< 3/4$. We suggest a method based on the higher criticism
that is powerful in the whole range $\alpha> 1/2$. This optimality
property is true for a variety of designs, including the classical
(balanced) multi-way designs and more modern ``$p>n$'' designs arising
in genetics and signal processing. In addition to the standard
fixed effects model, we establish similar results for a random
effects model where the nonzero coefficients of the regression
vector are normally distributed.
\end{abstract}

%

\end{frontmatter}

\section{Introduction}
\label{secintro}

\subsection{The analysis of variance}

Testing whether a subset of covariates have any linear relationship
with a quantitative response has been a staple of statistical analysis
since Fisher introduced the analysis of variance (ANOVA) in the 1920s
\cite{MR0346954}. Fisher developed ANOVA in the context of
agricultural trials and the test has since then been one of the
central tools in the statistical analysis of experiments
\cite{MR2552961}. As a consequence, there are countless situations in
which it is routinely used, in particular, in the analysis of clinical
trials \cite{MR2154988} or in that of cDNA microarray experiments
\cite{churchill2002fundamentals,kerr2000analysis,slonim2002patterns},
to name just two important areas of biostatistics.

To begin with, consider the simplest design known as the one-way
layout,
\[
y_{ij} = \mu+ \tau_j + z_{ij},
\]
where $y_{ij}$ is the $i$th observation in group $j$, $\tau_j$ is the
main effect for the $j$th treatment, and the $z_{ij}$'s are
measurement errors assumed to be i.i.d. zero-mean normal variables.
The goal is of course to determine whether there is any difference
between the treatments. Formally, assuming there are $p$ groups, the
testing problem is
\begin{eqnarray*}
&H_0\dvtx\tau_1 = \tau_2 =\cdots= \tau_p = 0,&\\
&H_1\dvtx\mbox{at least one } \tau_j \neq0.&
\end{eqnarray*}
The classical one-way analysis of variance is based on the
well-known $F$-test calculated by all statistical software packages. A
characteristic of ANOVA is that it tests for a \textit{global} null and
does not result in the identification of which $\tau_j$'s are nonzero.

Taking within-group averages reduces the model to
%
%
\begin{equation}
\label{eqbeta+z}
y_j = \beta_j + z_j,\qquad j = 1, \ldots, p,
\end{equation}
where $\beta_j = \mu+ \tau_j$ and the $z_j$'s are independent
zero-mean Gaussian variables. If we suppose that the grand mean has
been removed, so that the overall mean effect vanishes, that is, $\mu=
0$, then the testing problem becomes
%
%
\begin{eqnarray}
\label{eqtest}
&H_0\dvtx\beta_1 = \beta_2 =\cdots= \beta_p = 0,&\\
&H_1\dvtx\mbox{at least one } \beta_j \neq0.&\nonumber
\end{eqnarray}
In order to discuss the power of ANOVA in this setting, assume for
simplicity that the variances of the error terms in (\ref{eqbeta+z})
are known and identical, so that ANOVA reduces to a chi-square test
that rejects for large values of $\sum_j y_j^2$. As explained before,
this test does not identify which of the $\beta_j$'s are nonzero, but
it has great power in the sense that it maximizes the minimum power
against alternatives of the form $\{\bbeta\dvtx\sum_j \beta_j^2 \geq
B\}$ where $B > 0$. Such an appealing property may be shown via
invariance considerations; see \cite{MR941007} and
\cite{TSH}, Chapters 7 and 8.

\subsection{Multiple testing and sparse alternatives}

A different approach to the same testing problem is to test each
individual hypothesis $\beta_j = 0$ versus $\beta_j \neq0$, and
combine these tests by applying a Bonferroni-type correction. One way
to implement this idea is by computing the minimum $P$-value and
comparing it with a threshold adjusted to achieve a desired significance
level. When the variances of the $z_j$'s are identical, this is
equivalent to rejecting the null when
%
%
\begin{equation}
\label{eqmaxtest}
\operatorname{Max}(\by) = {\max_j} |y_j|
\end{equation}
exceeds a given threshold. From now on, we will refer to this
procedure as the Max test. Because ANOVA is such a well established
method, it might surprise the reader---but not the specialist---to
learn that there are situations where the Max test, though apparently
naive, outperforms ANOVA by a wide margin. Suppose indeed that $z_j
\sim\mathcal{N}(0,1)$ in (\ref{eqbeta+z}) and consider an alternative
of the form $\max_j |\beta_j| \geq A$ where $A > 0$. In this setting,
ANOVA requires $A$ to be at least as large as $p^{1/4}$ to provide
small error probabilities, whereas the Max test only requires $A$ to be
on the order of $(2\log p)^{1/2}$. When $p$ is large, the difference
is very substantial. Later in the paper, we shall prove that in an
asymptotic sense, the Max test maximizes the minimum power against
alternatives of this form. The key difference between these two
different classes of alternatives resides in the kind of
configurations of parameter values which make the likelihoods under
$H_0$ and $H_1$ very close. For the alternative $\{ \bbeta\dvtx\sum_j
\beta_j^2 \geq B\}$, the likelihood functions are hard to distinguish
when the entries of $\bbeta$ are of about the same size (in absolute
value). For the other, namely, $\{\bbeta\dvtx\max_j |\beta_j| \geq A\}$,
the likelihood functions are hard to distinguish when there is a
single nonzero coefficient equal to $\pm A$.

Multiple hypothesis testing with sparse alternatives is now
commonplace, in particular, in computational biology where the data is
high-dimensional and we typically expect that only a few of the many
measured variables actually contribute to the response---only a few
assayed treatments may have a positive effect. For instance, DNA
microarrays allow the monitoring of expression levels in cells for
thousands of genes simultaneously. An important question is to decide
whether some genes are differentially expressed, that is, whether or
not there are genes whose expression levels are associated with a
response such as the absence/presence of prostate cancer. A typical
setup is that the data for the $i$th individual consists of a response
or covariate $y_i$ (indicating whether this individual has a specific
disease or not) and a gene expression profile $y_{ji}$, $1 \le j \le
p$. A standard approach consists in computing, for each gene $j$, a
statistic $T_j$ for testing the null hypothesis of equal mean
expression levels and combining them with some multiple hypothesis
procedure \cite{Dudoit03,Efron00}. A possible and simple model in
this situation may assume $T_j \sim\mathcal{N}(0,1)$ under the null
while $T_j \sim\mathcal{N}(\beta_j,1)$ under the alternative. Hence, we
are in our sparse detection setup since one typically expects only a
few genes to be differentially expressed. Despite the form of
the alternative, ANOVA is still a popular method for testing the global
null in such problems \mbox{\cite{kerr2000analysis,slonim2002patterns}}.

\subsection{This paper}

Our exposition has thus far concerned simple designs, namely,
the one-way layout or sparse mean model. This paper, however, is
concerned with a much more general problem: we wish to decide whether
or not a response depends linearly upon a few covariates. We thus
consider the standard linear model
%
%
\begin{equation} \label{setting}
\by= \bX\bbeta+ \bz
\end{equation}
with an $n$-dimensional response $\by= (y_1,\ldots, y_n)$, a data
matrix $\bX\in\bbR^{n \times p}$ (assumed to have full rank) and a
noise vector, assumed to be i.i.d. standard normal. The decision
problem (\ref{eqtest}) is whether all the $\beta_i$'s are zero or
not. We briefly pause to remark that statistical practitioners are
familiar with the ANOVA derived $F$-statistic---also known as the
model adequacy test---that software packages routinely provide for
testing $H_0$. Our concern, however, is not at all model adequacy but
rather we view the test of the global null as a detection problem. In
plain English, we would like to know whether there is signal or
whether the data is just noise. A more general problem is to test
whether a subset of coordinates of $\bbeta$ are all zero or not, and,
as is well known, ANOVA is in this setup the most popular tool for
comparing nested models. We emphasize that our results also apply to
such general model comparisons, as we shall see later.

There are many applications of high-dimensional setups in which a
response may depend upon only a few covariates. We give a few examples
in the life sciences and in engineering; there are, of course, many
others:
\begin{itemize}
\item\textit{Genetics.} A single nucleotide polymorphism (SNP) is a
form of DNA variation that occurs when at a single position in the
genome, multiple (typically two) different nucleotides are found
with positive frequency in the population of reference. One then
collects information about allele counts at polymorphic locations.
Almost all common SNPs have only two alleles so that one records a
variable $x_{ij}$ on individual $i$ taking values in $\{0, 1, 2\}$
depending upon how many copies of, say, the rare allele one
individual has at location $j$. One also records a quantitative
trait $y_i$. Then the problem is to decide whether or not this
quantitative trait has a genetic background. In order to scan the
entire genome for a signal, one needs to screen between 300,000 and
1,000,000 SNPs. However, if the trait being measured has a genetic
background, it will be typically regulated by a small number of
genes. In this example, $n$ is typically in the thousands while $p$
is in the hundreds of thousands. The standard approach is to test
each hypothesis $H_j\dvtx\beta_j \neq0$ by using a statistic
depending on the least-squares estimate $\hat{\beta}_j$ obtained by
fitting the simple linear regression model
%
%
\begin{equation}
\label{eqslr}
y_i = \hat{\beta}_0 + \hat{\beta}_j x_{ij} + r_{ij}.
\end{equation}
The global null is then tested by adjusting the significance level to
account for the multiple comparisons, effectively implementing a Max
test; see \cite{Chiararef1,mccarthy2008genome}, for example.

\item\textit{Communications.} A multi-user detection problem typically
assumes a linear model of the form (\ref{setting}), where the $j$th
column of $\bX$, denoted $\vct{x}_j$, is the channel impulse response
for user $j$ so that
the received signal from the $j$th user is $\beta_j \vct{x}_j$ (we have
$\beta_j = 0$ in case user $j$ is not sending any message).
Note that the mixing matrix $\bX$ is often modeled as random with
i.i.d. entries.
In a strong noise environment, we might be interested in knowing whether
information is being transmitted (some $\beta_j$'s are not zero) or
not. In some applications, it is reasonable to assume that only a
few users are transmitting information at any given time.
Standard methods include the matched filter detector, which
corresponds to the Max test applied to $\bX^T \by$, and linear
detectors, which correspond to variations of the ANOVA
$F$-test \cite{honig2009advances}.

\item\textit{Signal detection.} The most basic problem in signal
processing concerns the detection of a signal $S(t)$ from the data
$y(t) = S(t) + z(t)$ where $z(t)$ is white noise. When the signal is
nonparametric, a popular approach consists in modeling $S(t)$ as a
(nearly) sparse superposition of waveforms taken from a dictionary
$\bX$, which leads to our linear model (\ref{setting}) (the columns
of $\bX$ are elements from this dictionary). For instance, to detect
a multi-tone signal, one would employ a dictionary of sinusoids; to
detect a superposition of radar pulses, one would employ a
time-frequency dictionary \cite{Mallat93,MallatBook}; and to detect
oscillatory signals, one would employ a dictionary of chirping
signals. In most cases, these dictionaries are massively
overcomplete so that we have more candidate waveforms than the
number of samples, that is, $p > n$. Sparse signal detection problems
abound, for example the detection of cracks in materials
\cite{Zhang2000961}, of hydrocarbon from seismic data
\cite{castagna120} and of tumors in medical imaging
\cite{breast-tumor}.

\item\textit{Compressive sensing.} The sparse detection model may also
arise in the area of compressive sensing
\cite{CRT,OptimalRecovery,Donoho-CS}, a novel theory which asserts
that it is possible to accurately recover a (nearly) sparse
signal---and by extension, a signal that happens to be sparse in some
fixed basis or dictionary---from the knowledge of only a few of its
random projections. In this context, the $n \times p$ matrix $\bX$
with $n \ll p$ may be a random projection such as a partial Fourier
matrix or a matrix with i.i.d. entries. Before reconstructing the
signal, we might be interested in testing whether there is any
signal at all in the first place.
\end{itemize}

All these examples motivate the study of two classes of sparse
alternatives:
\begin{longlist}[(1)]
\item[(1)] \textit{Sparse fixed effects model} (\textit{SFEM}). Under the alternative,
the regression vector $\bbeta$ has at least $S$ nonzero coefficients
exceeding $A$ in absolute value.
\item[(2)] \textit{Sparse random effects model} (\textit{SREM}). Under the
alternative, the regression vector $\bbeta$ has at least $S$ nonzero
coefficients assumed to be i.i.d. normal with zero mean and variance
$\tau^2$.
\end{longlist}
In both models, we set $S = p^{1-\alpha}$, where $\alpha\in(0,1)$ is
the sparsity exponent. Our purpose is to study the performance of
various test statistics for detecting such
alternatives.\setcounter{footnote}{1}\footnote{We
will sometimes put a prior on the support of $\bbeta$ and on the
signs of its nonzero entries in SFEM.}

\subsection{Prior work}
\label{secpriorWork}
To introduce our results and those of others, we need to recall a few
familiar concepts from statistical decision theory. From now on,
$\Omega$ denotes a set of alternatives, namely, a subset of $\bbR^p
\setminus\{0\}$ and $\pi$ is a prior on $\Omega$. The Bayes risk of a
test $T = T(\bX, \by)$ for testing $\bbeta= {\mathbf0}$ versus $\bbeta
\sim\pi$ when $H_0$ and $H_1$ occur with the same probability is
defined as the sum of its probability of type I error (false alarm)
and its average probability of type II error (missed
detection). Mathematically,
%
%
\begin{equation} \label{risk}
\operatorname{Risk}_\pi(T) := \mathbb{P}_{\mathbf0}(T = 1) + \pi[\mathbb
{P}_\bbeta(T =
0)],
\end{equation}
where $\mathbb{P}_\bbeta$ is the probability distribution of $\by$
given by
the model (\ref{setting}) and $\pi[\cdot]$ is the expectation with
respect to the prior $\pi$. If we consider the linear model in the
limit of large dimensions, that is, $p \to\infty$ and $n = n(p) \to
\infty$,
and a sequence of priors $\{\pi_p\}$, then we say that a sequence of
tests $\{T_{n,p}\}$ is
asymptotically \textit{powerful} if $\lim_{p \to\infty}
\operatorname{Risk}_{\pi_p}(T_{n,p}) = 0$. We say that it is
asymptotically \textit{powerless} if
$\liminf_{p \to\infty} \operatorname{Risk}_{\pi_p}(T_{n,p}) \geq1$.
When no
prior is specified, the risk is understood as the worst-case risk
defined as
\[
\operatorname{Risk}(T) := \mathbb{P}_{\mathbf0}(T = 1) + \max_{\bbeta\in
\Omega}
\mathbb{P}_\bbeta(T = 0).
\]

With our modeling assumptions, ANOVA for testing $\bbeta= {\mathbf0}$
versus $\bbeta\neq{\mathbf0}$ reduces to the chi-square test that
rejects for large values of $\|\mtx{P} \vct{y}\|^2$, where $\mtx{P}$
is the orthogonal projection onto the range of $\mtx{X}$. Since under
the alternative, $\|\mtx{P} \vct{y}\|^2$ has the chi-square
distribution with $\min(n,p)$ degrees of freedom and noncentrality
parameter $\|\bX\bbeta\|^2$, a simple argument shows that ANOVA is
asymptotically powerless when
%
%
\begin{equation} \label{anova}
\|\bX\bbeta\|^2/\sqrt{\min(n,p)} \to0,
\end{equation}
and asymptotically powerful if the same quantity tends to infinity.
This is congruent with the performance of ANOVA in a standard one-way
layout; see \cite{MR2062823}, who obtain the weak limit of the ANOVA
$F$-ratio under various settings.

Consider the sparse fixed effects alternative now. We prove that ANOVA
is still essentially optimal under mild levels of sparsity
corresponding to $\alpha\in[0,1/2]$ but not under strong sparsity
where $\alpha\in(1/2,1]$. In the sparse mean model
(\ref{eqbeta+z}) where $\bX$ is the identity, ANOVA is suboptimal,
requiring $A$ to grow as a power of $p$; this is simply because
(\ref{anova}) becomes $A^2 S/\sqrt{p} \to0$ when all the nonzero
coefficients are equal to $A$ in absolute value. In contrast, the Max
test is asymptotically powerful when $A$ is on the order of
$\sqrt{\log p}$ but is only optimal under very strong sparsity,
namely, for $\alpha\in[3/4,1]$. It is possible to improve on the Max
test in the range $\alpha\in(1/2, 3/4)$ and we now review the
literature which only concerns the sparse mean model, $\bX=
\bI_p$. Set
%
%
\begin{equation}
\label{eqthreshold}
\rho^*(\alpha) = \cases{
\alpha-1/2, &\quad$1/2 < \alpha< 3/4$, \vspace*{2pt}\cr
\bigl(1 -\sqrt{1 -\alpha}\bigr)^2, &\quad$3/4 \leq\alpha< 1$.}
\end{equation}
Then Ingster \cite{Ingster99} showed that if $A = \sqrt{2 r \log p}$
with $r < \rho^*(\alpha)$ fixed as $p \to\infty$, then all sequences
of tests are asymptotically powerless. In the other direction, he
showed that there is an asymptotically powerful sequence of tests if
$r > \rho^*(\alpha)$. See also the work of Jin \cite{jinPhD}. Donoho
and Jin
\cite{dj04} analyzed a number of testing procedures in this setting,
and, in particular, the higher criticism of Tukey which rejects for
large values of
\[
\mathrm{H C}^*(\by) = \sup_{t > 0} \frac{\# \{i\dvtx|y_i| > t\} -
2p\bar{\Phi}(t)}{\sqrt{2 p\bar{\Phi}(t) (1 -2 \bar{\Phi}(t))}},
\]
where $\bar{\Phi}$ denotes the survival function of a standard normal
random variable. They showed that the higher criticism is powerful
within the detection region established by Ingster. Hall and Jin
\cite{hj08,hj09} have recently explored the case where the noise may
be correlated, that is, $\bz\sim\cN({\mathbf0}, \mtx{V})$ and the
covariance matrix $\mtx{V}$ is known and has full rank. Letting
$\mtx{V} = \mtx{L} \mtx{L}^T$ be a Cholesky factorization of the
covariance matrix, one can whiten the noise in $\by= \bbeta+ \bz$ by
multiplying both sides by $\mtx{L}^{-1}$, which yields $\tilde{\by} =
\mtx{L}^{-1}\bbeta+ \tilde{\bz}$; $\tilde{\bz}$ is now white noise,
and this is a special case of the linear model (\ref{setting}). When
the design matrix is triangular with coefficients decaying
polynomially fast away from the diagonal, \cite{hj09} proves that the
detection threshold remains unchanged, and that a form of higher criticism still achieves asymptotic optimality.

There are few other theoretical results in the literature, among which
\cite{MR2278336} develops a locally most powerful (score) test in a
setting similar to SREM; here, ``locally'' means that this property only
holds for values of $\tau$ sufficiently close to zero. The authors do
not provide any minimal value of $\tau$ that would guarantee the
optimality of their method. However, since their score test resembles
the ANOVA $F$-test, we suggest that it is only optimal for very small
values of $\tau$ corresponding to mild levels of sparsity, that
is, $\alpha< 1/2$.


Since the submission of our paper, a manuscript by Ingster, Tsybakov
and Verzelen \cite{ingster2010detection}, also considering the
detection of a sparse vector in the linear regression model, has become
publicly available. We comment on differences in Section \ref{secsparse}.

In the signal processing literature, a number of applied papers
consider the problem of detecting a signal expressed as a linear
combination in a dictionary
\cite{castagna120,Zhang2000961,MR1956096}. However, the
extraction of the salient signal is often the end goal of real signal
processing applications so that research has focused on estimation
rather than pure detection. As a consequence, one finds a literature
entirely focused on estimation rather than on testing whether the data
is just white noise or not. Examples of pure detection papers include
\cite{sparse-detection,CS-detection,meng-sparse}. In
\cite{sparse-detection}, the authors consider detection by matched
filtering, which corresponds to the Max test, and perform simulations
to assess its power. The authors in \cite{CS-detection} assume that
$\bbeta$ is approximately known and examine the performance of the
corresponding matched filter.
Finally, the paper \cite{meng-sparse} proposes a Bayesian approach for the
detection of sparse signals in a sensor network for which the design
matrix is assumed to have some polynomial decay in terms of the
distance between sensors.

\subsection{Our contributions}
\label{secpeak}

We show that if the predictor variables are not too correlated, there
is a sharp detection threshold in the sense that no test is
essentially better than a coin toss when the signal strength is below
this threshold, and that there are statistics which are asymptotically
powerful when the signal strength is above this threshold. This
threshold is the same as that one gets for the sparse mean
problem. Therefore, this work extends the earlier results and
methodologies cited above \cite{Ingster99,jinPhD,dj04,hj08,hj09}, and
is applicable to the modern high-dimensional situation where the
number of predictors may greatly exceed the number of observations.

A simple condition under which our results hold is a low-coherence
assumption.\footnote{Although we are primarily interested in the
modern $p > n$
setup, our results apply \textit{regardless} of the values of $p$ and
$n$.}
Let $\vct{x}_1,\ldots, \vct{x}_p$ be the column vectors of $\bX$,
assumed to
be normalized; this assumption is merely for convenience since it
simplifies the exposition, and is not essential. Then if a large
majority of all pairs of predictors have correlation less than
$\gamma$ with $\gamma= O(p^{-1/2 + \eps})$ for each $\eps> 0$ (the
real condition is weaker), then the results for the sparse mean model
(\ref{eqbeta+z}) apply almost unchanged. Interestingly, this is true
even when the ratio between the number of observations and the number
of variables is negligible, that is, $n/p \to0$. In particular, $A =
\sqrt{2\rho^*(\alpha) \log p}$ is the sharp detection threshold for
SFEM (sparse fixed effects model). Moreover, applying the higher criticism, not to the values of $\by$, but to those of $\bX^T \by$ is
asymptotically powerful as soon as the nonzero entries of $\bbeta$ are
above this threshold; this is true for all $\alpha\in(1/2, 1]$. In
contrast, the Max test applied to $\bX^T \by$ is only optimal in the
region $\alpha\in[3/4, 1]$. We derive the sharp threshold for SREM
as well, which is at $\tau= \sqrt{\alpha/(1-\alpha)}$. We show that
the Max tests and the higher criticism are essentially optimal in this
setting as well for all $\alpha\in(1/2, 1]$, that is, they are both
asymptotically powerful as soon as the signal-to-noise ratio permits.

Before continuing, it may be a good idea to give a few examples of
designs obeying the low-coherence assumption (weak correlations
between most of the predictor variables) since it plays an important
role in our analysis:
\begin{itemize}
\item\textit{Orthogonal designs.} This is the situation where the
columns of $\bX$ are orthogonal so that $\bX^T\bX$ is the $p \times
p$ identity matrix (necessarily, $p \leq n$). Here the coherence is
of course the lowest since $\gamma(\bX) = 0$.

\item\textit{Balanced, one-way designs.} As in a clinical trial
comparing $p$ treatments, assume a balanced, one-way design with $k$
replicates per treatment group and with the grand mean already
removed. This corresponds to the linear model (\ref{setting}) with
$n = p k$ and, since we assume the predictors to have norm~$1$,
%
%
\begin{equation} \label{one-way}
\bX= \frac{1}{\sqrt{k}} \left[\matrix{
{\mathbf1} & {\mathbf0} & \cdots& {\mathbf0} \cr
{\mathbf0} & {\mathbf1} & \cdots& {\mathbf0} \cr
\vdots& \vdots& \vdots& \vdots\cr
{\mathbf0} & {\mathbf0} & \cdots& {\mathbf1}}\right]
\in\bbR^{n \times p},
\end{equation}
where each vector in this block representation is $k$-dimensional.
This is in fact an example of orthogonal design. Note\vspace*{1pt} that our
results apply even under the standard constraint ${\mathbf1}^T \bbeta=
0$.

\item\textit{Concatenation of orthonormal bases.} Suppose that $p = n
k$ and that $\bX$ is the concatenation of $k$ orthonormal bases in
$\bbR^n$ jointly used as to provide an efficient signal
representation. Then our result applies provided that $k =
O(n^{\eps}), \forall\eps> 0$ and that our bases are mutually
incoherent so that $\gamma$ is sufficiently small (for examples of
incoherent bases see, e.g., \cite{MR1872845}).

\item\textit{Random designs.} As in some compressive sensing and
communications applications, assume that $\bX$ has i.i.d. normal
entries\footnote{This is a frequently discussed channel model in
communications.} with columns subsequently normalized (the column
vectors are sampled independently and uniformly at random on the
unit sphere). Such a design is close to orthogonal since $\gamma
\leq\sqrt{5 (\log p)/n}$ with high probability. This fact follows
from a well-known concentration inequality for the uniform
distribution on the sphere \cite{MR1849347}.
The exact same bound applies if the entries of $\bX$ are instead
i.i.d. Rademacher random variables.
\end{itemize}

We return to the discussion of our statistics and note that the higher criticism and the Max test applied to $\bX^T \by$ are exceedingly
simple methods with a straightforward implementation running in
$O(np)$ flops. This brings us to two important points:
\begin{longlist}[(1)]
\item[(1)] In the classical sparse mean model, Bonferroni-type multiple
testing (the Max test) is not optimal when the sparsity level is
moderately strong, that is, when $1/2 < \alpha< 3/4$ \cite{dj04}. This
has direct implications in the fields of genetics and genomics where
this is the prevalent method. The same is true in our more general
model and it implies, for example, that the matched filter detector
in wireless multi-user detection is suboptimal in the same sparsity
regime.

We elaborate on this point because this carries an important
message. When the sparsity level is moderately strong, the higher criticism method we propose is powerful in situations where the
signal amplitude is so weak that the Max test is powerless.
\textit{This says that one can detect a linear relationship between a
response $\by$ and a few covariates even though those covariates
that are most correlated with $\by$ are not even in the model.}
Put differently, if we assign a $P$-value to each hypothesis
$\beta_j = 0$ (computed from a simple linear regression as discussed
earlier), then \textit{the case against the null is not in the tail of
these $P$-values but in the bulk}, that is, the smallest
$P$-values may not carry any information about the presence of a
signal. In the situation we describe, the smallest $P$-values most
often correspond to true null hypotheses, sometimes in such a way
that the false discovery rate (FDR) cannot be controlled at any
level below 1; and yet, the higher criticism has full power.

\item[(2)] Though we developed the idea independently, the higher criticism
applied to $\bX^T \by$ is similar to the innovated higher criticism
of Hall and Jin \cite{hj09}, which is specifically designed for time
series. Not surprisingly, our results and arguments bear some
resemblance with those of Hall and Jin \cite{hj09}. We have already
explained how their results apply when the design matrix is
triangular (and, in particular, square) and has sufficiently rapidly
decaying coefficients away from the diagonal. Our results go much
further in the sense that (1) they include designs that are far from
being triangular or even square, and (2) they include designs with
coefficients that do not necessarily follow any ordered decay
pattern. On the technical side, Hall and Jin astutely reduce
matters to the case where the design matrix is banded, which greatly
simplifies the analysis. In the general linear model, it is not
clear how a similar reduction would operate especially when
$n < p$---at the very least, we do not see a way---and one must deal with
more intricate dependencies in the noise term $\bX^T \bz$.
\end{longlist}

As we have remarked earlier, we have discussed testing the global null
$\bbeta= {\mathbf0}$, whereas some settings obviously involve nuisance
parameters as in the comparison of nested models. Examples of nuisance
parameters include the grand mean in a balanced, one-way design or,
more generally, the main effects or lower-order interactions in a
multi-way layout. In signal processing, the nuisance term may
represent clutter as opposed to noise. In general, we have
\[
\by= \bX^{(0)} \bbeta^{(0)} + \bX^{(1)} \bbeta^{(1)} + \bz,
\]
where $\bbeta^{(0)}$ is the vector of nuisance parameters, and
$\bbeta^{(1)}$ the vector we wish to test. Our results concerning the
performance of ANOVA, the higher criticism or the Max test apply
provided that the column spaces of $\bX^{(0)}$ and $\bX^{(1)}$ be
sufficiently far apart. This occurs in lots of applications of
interest. In the case of the balanced, multi-way design, these spaces
are actually orthogonal. In signal processing, these spaces will also
be orthogonal if the column space of $\bX^{(0)}$ spans the
low-frequencies while we wish to detect the presence of a
high-frequency signal. The general mechanism which allows us to
automatically apply our results is to simply assume that $\mtx{P}_0
\bX^{(1)}$, where $\mtx{P}_0$ is the orthogonal\vspace*{1pt} projector with the
range of $\bX^{(0)}$ as null space, obeys the conditions we have for
$\bX$.

\subsection{Organization of the paper}

The paper is organized as follows. In Section~\ref{secortho} we consider
orthogonal designs and state results for the classical setting where
no sparsity assumption is made on the regression vector $\bbeta$, and
the setting where $\bbeta$ is mildly sparse. In Section \ref
{secsparse} we study
designs in which \textit{most} pairs of predictor variables are only
weakly correlated; this part contains our main results. In
Section \ref{secfull} we focus on some examples of designs with full
correlation structure, in particular, multi-way layouts with embedded
constraints. Section \ref{secnumerics} complements our study with some
numerical experiments, and we close the paper with a short discussion,
namely, Section \ref{secdiscussion}. Finally, the proofs are gathered
in a
supplementary file \cite{unstructured-suppl}.

\subsection{Notation}

We provide a brief summary of the notation used in the paper.
Set $[p] = \{1,\ldots, p\}$ and for a subset $\cJ\subset[p]$, let
$|\cJ|$ be its cardinality. Bold upper (resp., lower) case letters
denote matrices (resp., vectors), and the same letter not bold
represents its coefficients, for example, $a_j$ denotes the $j$th entry of
$\vct{a}$. For an $n \times p$ matrix $\bA$ with column vectors $\vct
{a}_1,\ldots, \vct{a}_p$, and a subset $\cJ\subset[p]$, $\bA_\cJ$
denotes the
$n$-by-$|\cJ|$ matrix with column vectors $\vct{a}_j, j \in\cJ$.
Likewise, $\vct{a}_\cJ$ denotes the vector $(a_j, j \in\cJ)$. The
Euclidean norm of a vector is $\|\vct{a}\|$ and the sup-norm
$\|\vct{a}\|_\infty$. For a matrix $\bA= (a_{ij})$, $\|\bA\|_\infty=
\sup_{i,j} |a_{ij}|$, and this needs to be distinguished from
$\|\bA\|_{\infty, \infty}$, which is the operator norm induced by the
sup norm, $\|\bA\|_{\infty, \infty} = \sup_{\|\vct{x}\|_\infty\le1}
\|\bA\vct{x}\|_\infty$. The Frobenius (Euclidean) norm of $\bA$ is
$\|\bA\|_F$. $\Phi$ (resp., $\phi$) denotes the cumulative distribution
(resp., density) function of a standard normal random variable, and
$\bar{\Phi}$ its survival function. For brevity, we say that $\bbeta$
is $S$-sparse if $\bbeta$ has exactly $S$ nonzero
coefficients. Finally, we say that a random variable $X \sim F_X$ is
stochastically smaller than $Y \sim F_Y$, denoted $X
\leq^{\mathrm{sto}}
Y$, if $F_X(t) \geq F_Y(t)$ for all scalar $t$.

\section{Orthogonal designs}
\label{secortho}

This section introduces some results for the orthogonal design in
which the columns of $\bX$ are orthonormal, that is, $\bX^T \bX=
\bI_p$. While from the analysis viewpoint there is little difference
with the case where $\bX$ is the identity matrix, this is of course a
special case of our general results, and this section may also serve
as a little warm-up. Our first result, which is a special case of
Proposition \ref{prpmoderate-lb}, determines the range of sparse
alternatives for
which ANOVA is essentially optimal.
\begin{prp} \label{prpmoderate}
Suppose $\bX$ is orthogonal and let the number of nonzero coefficients
be $S = p^{1 - \alpha}$ with $\alpha\in[0,1/2]$. In SFEM
(resp., SREM), all sequences of tests are asymptotically powerless if
$A^2 S/p^{1/2} \to0$ (resp.,\break $\tau^2 S/p^{1/2} \to0$).
\end{prp}

Returning to our earlier discussion, it follows from (\ref{anova}) and
the lower bound $\|\bX\beta\|^2 = \|\beta\|^2 \ge A^2 S$ that
ANOVA\vspace*{1pt}
has full asymptotic power whenever $A^2 S/p^{1/2} \to
\infty$. Therefore, comparing this with the content of Proposition
\ref{prpmoderate} reveals that ANOVA is essentially optimal in the
moderately sparse range corresponding to $\alpha\in[0,1/2]$.

The second result of this section is that under
an $n \times p$ orthogonal design, the detection threshold is the same
as if $\bX$ were the identity. We need a little bit of notation to
develop our results. As in \cite{dj04}, define
\[
\rho_{\mathrm{Max}}(\alpha) = \bigl(1 - \sqrt{1-\alpha}\bigr)^2,
\]
and observe that with $\rho^*(\alpha)$ as in (\ref{eqthreshold}),
\[
\cases{
\rho^*(\alpha) < \rho_{\mathrm{Max}}(\alpha), &\quad$1/2 \le\alpha<
3/4$,\cr
\rho^*(\alpha) = \rho_{\mathrm{Max}}(\alpha), &\quad$3/4 \le\alpha\le1$.}
\]
We will also set a detection threshold for SREM defined by
%
%
\begin{equation}
\label{eqthreshrd}
\rho^*_{\mathrm{rand}}(\alpha) = \sqrt{\alpha/(1 -\alpha)}.
\end{equation}
With these definitions, the following theorem compares the performance
of the higher criticism and the Max test.
\begin{theorem} \label{thmortho} Suppose $\bX$ is orthogonal and assume
the sparsity exponent obeys $\alpha\in(1/2,1]$.
\begin{longlist}[(1)]
\item[(1)] In SFEM,\vspace*{1pt} all sequences of tests are asymptotically powerless
if $A = \sqrt{2 r \log p}$ with $r < \rho^*(\alpha)$. Conversely,
the higher criticism applied to $|\vct{x}_1^T\by|,\ldots,\break |\vct{x}_p^T
\by|$ is asymptotically powerful if $r > \rho^*(\alpha)$. Also,
the Max test is asymptotically powerful if $r > \rho_{\mathrm
{Max}}(\alpha)$ and powerless if $r < \rho_{\mathrm{Max}}(\alpha)$.

\item[(2)] In SREM, all sequences of tests are asymptotically powerless
if $\tau< \rho^*_{\mathrm{rand}}(\alpha)$. Conversely,\vspace*{1pt} both the higher criticism and the Max test applied to $|\vct{x}_1^T\by|,\ldots,
|\vct{x}_p^T \by|$ are asymptotically powerful if $\tau> \rho^*_{\mathrm
{rand}}(\alpha)$.
\end{longlist}
In the upper bounds, $r$ and $\tau$ are fixed while $p \to\infty$.
\end{theorem}

To be absolutely clear, the statements for SFEM may be understood
either in the worst-case risk sense or under the uniform prior on the
set of $S$-sparse vectors with nonzero coefficients equal to $\pm
A$. For SREM, the prior simply selects the support of $\bbeta$
uniformly at random.
After multiplying the observation by $\bX^T$,
matters are reduced to the case of the identity design for which the
performance of the higher criticism and the Max test have been
established in SFEM \cite{dj04}. The result for the sparse random
model is new and appears in more generality in Theorem \ref{thmmax}.

To conclude, the situation concerning orthogonal designs is very
clear. In SFEM, for instance, if the sparsity level is such that
$\alpha\le1/2$, then ANOVA is asymptotically optimal whereas the
higher criticism is optimal if $\alpha> 1/2$. In contrast, the Max
test is only optimal in the range $\alpha\geq3/4$.

\section{Weakly correlated designs}
\label{secsparse}

We begin by introducing a model of design matrices in which most of
the variables are only weakly correlated. Our model depends upon two
parameters, and we say that a $p \times p$ correlation matrix $\vct{C}$
belongs to the class $\cS_p(\gamma, \Delta)$ if and only if it
obeys the following two properties:
\begin{itemize}
\item\textit{Strong correlation property.} This requires that for all $j
\neq k$,
\[
|c_{jk}| \le1-(\log p)^{-1}.
\]
That is, \textit{all} the correlations are bounded above by $1- (\log
p)^{-1}$. In the limit of large $p$, this is not an assumption and we
will later explain how one can relax this even further.

\item\textit{Weak correlation property.} This is the main assumption and
this requires that for all $j$,
\[
\bigl|\{k \dvtx|c_{jk}| > \gamma\}\bigr| \leq\Delta.
\]
Note that for $\gamma\le1$, $\Delta\geq1$ since $c_{jj} = 1$. Fix
a variable $\vct{x}_j$. Then at most $\Delta- 1$ other variables have a
correlation exceeding $\gamma$ with $\vct{x}_j$.
\end{itemize}
Our only real condition caps the number of variables that can
have a correlation with any other above a threshold $\gamma$. An
orthogonal design belongs to $\cS_p(0,1)$ since all the correlations
vanish. With high probability, the Gaussian and Rademacher designs
described earlier belong to $\cS_p(\gamma,1)$ with $\gamma= \sqrt{5
(\log p)/n}$.

\subsection{Lower bound on the detectability threshold}

The main result of this paper is that if the predictor variables are
not highly correlated, meaning that the quantities $\gamma$ and
$\Delta$ above are sufficiently small, then there are computable
detection thresholds for our sparse alternatives that are very similar
or identical to those available for orthogonal designs.

We begin by studying lower bounds and for SFEM, these may be
understood either in a worst-case sense or under the prior where
$\bbeta$ is uniformly distributed among all $S$-sparse vectors with
nonzero coefficients equal to $\pm A$. For SREM, these hold under a
prior generating the support uniformly at random. We first consider
mildly sparse alternatives.
\begin{prp} \label{prpmoderate-lb} 
Suppose that $\bX^T \bX\in\cS_p(\gamma,
1)$ and let $S = p^{1 - \alpha}$ with $\alpha\in[0,1/2]$.
In SFEM (resp., SREM), all sequences of tests are asymptotically powerless
if $A^2 S (p^{-1/2} + \gamma\log p) \to0$ [resp., $\tau^2 S (p^{-1/2} +
\gamma) \to0$].
\end{prp}

In order to interpret this proposition, we note that $\gamma$ will
usually be at least as large as $n^{-1/2}$, as shown just below.

In Proposition \ref{prpmoderate-lb} we have required that $\Delta= 1$
in order to
derive sharp results. Moving now to sparser alternatives, we allow for
$\Delta$ to increase with $p$, although very slowly, while the condition
on $\gamma$ remains essentially the same.
\begin{theorem} \label{thmSp-lb} Assume the sparsity exponent obeys
$\alpha\in(1/2,1]$, and suppose that $\bX^T \bX\in\cS_p(\gamma,
\Delta)$ with the following parameter asymptotics: (1) $\Delta=
O(p^{\eps})$, for all $\eps> 0$, and (2) $\gamma p^{1-\alpha} (\log
p)^4 \to0$.
In SFEM (resp., SREM), all sequences of tests are asymptotically powerless
if $A = \sqrt{2 r \log p}$ with $r < \rho^*(\alpha)$ [resp., $\tau<
\rho^*_{\mathrm{rand}}(\alpha)$].
\end{theorem}

The result is essentially the same in the case of a balanced,
multi-way design with the usual linear constraints. We comment on
this point at the end of the proof of Theorem \ref{thmSp-lb}.

The reader may be surprised to see that the number $n$ of observations
does not explicitly appear in the above lower bounds. The sample size
appears implicitly, however, since it must be large enough for the
class $\cS_p(\gamma,\Delta)$ to be nonempty. Assume $\Delta= 1$, for
instance, and that $p \geq n$.
Then by the lower bound \cite{MR1984549}, equation (12), we have
%
%
\begin{equation} \label{gamma-lb}
\gamma\geq\sqrt{(p-n)/(n p)}.
\end{equation}
For instance, $\gamma\geq1/\sqrt{2 n}$ if $p \geq2 n$.

As a technical aside, we remark that the lower bounds hold under the
strong correlation assumption
\[
|c_{jk}| \le1-\delta
\]
for any $\delta< 1$, provided that $\gamma\delta^{-2} p^{1-\alpha}
(\log p)^{3/2} \to0$. We shall prove this more general statement, and
the theorem is thus a special case corresponding to $\delta= (\log
p)^{-1}$.

We pause to compare with the results of the recent paper \cite
{ingster2010detection}. The lower bounds in \cite
{ingster2010detection} are the same as ours (for SFEM) except that they
impose slightly weaker conditions on $\gamma$. In Proposition \ref
{prpmoderate-lb}, their condition is $A^2 S (p^{-1/2} + \gamma) \to0$,
and in Theorem \ref{thmSp-lb}, their condition is $\gamma p^{1-\alpha
} \log p
\to0$.

\subsection{Upper bound on the detectability threshold}

We now turn to upper bounds and, unless stated otherwise, these assume
the following models:
\begin{itemize}
\item For SFEM, we assume that $\bbeta$ has a support generated
uniformly at random and that its nonzero coefficients have random
signs.
\item For SREM, we assume that $\bbeta$ has a support generated
uniformly at random.
\end{itemize}
We require that the support of $\bbeta$ be generated uniformly at
random and, in SFEM, that the signs of its coefficients be also random
to rule out situations where cancellations occur, making the signal
strength potentially too small (and possibly vanish) to allow for
reliable detection.

We begin by studying the performance of ANOVA when the alternative is
not that sparse. We state our result for $\Delta= 1$ in accordance
with the lower bound (Proposition \ref{prpmoderate-lb}), although the
result holds
when $\Delta$ obeys $\Delta= O(p^\eps)$ for all $\eps> 0$.
\begin{prp} \label{prpmoderate-ub} 
Assume that $\bX^T \bX\in\cS_p(\gamma, 1)$ and let $S = p^{1 -
\alpha}$.
\begin{itemize}
\item Assume $\gamma\log p \to0$. Then, in SFEM, ANOVA is
asymptotically powerful (resp., powerless) when $A^2 S/\sqrt{\min(n,
p)} \rightarrow\infty$ (resp., $\to0$).
%
\item Assume $\gamma\to0$. Then, in SREM, ANOVA is asymptotically
powerful (resp., powerless) when $\tau^2 S/\sqrt{\min(n, p)}
\rightarrow\infty$ (resp., $\to0$).

\end{itemize}
Note that this holds for all values of $\alpha$.
\end{prp}

For example, consider an $n \times p$ Gaussian design with $p > n$.
For this design $\gamma\asymp\sqrt{(\log p)/n}$ (in probability).
Hence, assuming $(\log p)^{3/2}/\sqrt{n} \to0$, Proposition \ref
{prpmoderate-ub}
says that, in SFEM, the ANOVA test is powerful when $A^2 S/\sqrt{n}
\rightarrow\infty$. We contrast this with Proposition \ref{prpmoderate-lb},
which says that, in the same context and assuming that $\alpha\in
[0,1/2]$, all methods are powerless when $A^2 S (\log
p)^{3/2}/\sqrt{n} \rightarrow0$. Hence, in this moderately sparse
setting where $\alpha\in[0,1/2]$, if one ignores the $(\log
p)^{3/2}$ factor (we do not know whether Proposition \ref
{prpmoderate-lb} is
tight), then one can say that ANOVA achieves the optimal detection
boundary.
However, as we will see in Theorems \ref{thmSp-ub}, \ref{thmSp-ub2}
and \ref{thmmax}, ANOVA is far from optimal in the strongly sparse
case when $\alpha> 1/2$.

Compared with Proposition \ref{prpmoderate-lb}, the condition on
$\gamma$ is
substantially weaker. More importantly, there appears to be a major
discrepancy when $n$ is negligible compared to $p$ because
$\sqrt{\min(n,p)}$ replaces $\sqrt{p}$. This is illusory, however, as
the lower bound on $\gamma$ displayed in (\ref{gamma-lb}) implies that
the condition on $A$ in Proposition~\ref{prpmoderate-lb} matches that of
Proposition \ref{prpmoderate-ub} up to a $\log p$ factor.

Turning to sparser alternatives, we apply the higher criticism to
$\bX^T \by$ and for $t > 0$, put
\[
H(t) = \frac{|\{j \dvtx|\vct{x}_j^T \by| > t\}| - 2p\bar{\Phi
}(t)}{\sqrt{2
p\bar{\Phi}(t)(1-2\bar{\Phi}(t))}}.
\]
The innovated higher criticism of Hall and Jin \cite{hj09} resembles
$\sup_{t>0} H(t) := \mathrm{H C}^*(\bX^T\by)$, the main difference being that
they apply a threshold to the entries of $\bX$ before multiplying by
$\bX^T$. Here, to facilitate the analysis, we search for the maximum
on a discrete grid and define
\[
H^*(s) = \max \bigl\{H(t) \dvtx t \in\bigl[s, \sqrt{5 \log p}\bigr] \cap\bbN\bigr\}.
\]

\begin{theorem} \label{thmSp-ub} Assume the sparsity exponent obeys
$\alpha\in(1/2,1]$ and suppose that $\bX^T \bX\in\cS_p(\gamma,
\Delta)$
with the following parameter asymptotics: (1) $\Delta= O(p^{\eps})$,
for all $\eps> 0$; (2) $\gamma^2 p^{1-\alpha} (\log p)^3 \to0$ and
(3) $\gamma^3 = O(p^{\eps+ 5 \alpha-4})$, for all $\eps> 0$.
\begin{itemize}
\item In SFEM, the test based on $H^*(\sqrt{2 r_\alpha\log p})$ with
$r_\alpha:= \min(1, 4 \rho^*(\alpha))$ is
asymptotically powerful against any alternative defined by $S =
p^{1-\alpha'}$ with $\alpha' \geq\alpha$ and $A = \sqrt{2 r \log
p}$ with $r > \rho^*(\alpha')$.
\item In SREM, the test based on $H^*(\sqrt{2 \log p})$ is
asymptotically powerful when $\tau> \rho^*_{\mathrm{rand}}(\alpha)$
regardless of $\alpha\in(1/2,1]$ and without condition (3).
\end{itemize}
\end{theorem}

In SREM, the conclusion is an immediate consequence of the behavior of
the Max test stated in Theorem \ref{thmmax} and we, therefore, omit
the proof.
Having said this, the remarks below apply to SFEM:
\begin{longlist}[(1)]
\item[(1)] The condition on $\gamma$ is weaker than the condition required
in Theorem \ref{thmSp-lb}, although the two conditions get ever closer
as $\alpha$ approaches $1/2$.

\item[(2)] The test based on $H^*(\sqrt{2 \log p})$ is asymptotically
powerful for all $\alpha\in[3/4, 1]$ (this test is closely related
to the Max test).


\item[(3)] Other discretizations in the definition of $H^*$ would yield the
same result. In fact, we believe the result holds without any
discretization, but we were not able to establish this in general.
However, suppose that $p = k n$ and that $\bX$ is the concatenation
of $k$ orthonormal bases. If $k = O(n^{\eps})$, for all $\eps> 0$,
the result holds without any discretization, meaning that rejecting
for large values of $\sup_{t>0} H(t)$ is asymptotically powerful
under the same conditions. This comes from leveraging the behavior
(under the null) of the higher criticism---detailed in
\cite{dj04}---for each basis.
\end{longlist}

While the above theorem gives relatively weak requirements on $\gamma
$, it is not fully adaptive. In particular, in SFEM, one requires
knowledge of $\alpha$ to set the search grid for the statistic $H^*$.
Under a stronger condition on $\gamma$, we have the following fully
adaptive result for $\alpha\in(1/2,1]$.
\begin{theorem} \label{thmSp-ub2} Assume the sparsity exponent obeys
$\alpha\in(1/2,1]$ and suppose that $\bX^T \bX\in\cS_p(\gamma,
\Delta)$
with the following parameter asymptotics: (1) $\Delta= O(p^{\eps})$,
for all $\eps> 0$; (2) $\gamma= O(p^{-1/2+\eps})$, for all $\eps>
0$. Then in SFEM, the test based on $H^*(1)$ is asymptotically
powerful whenever $r > \rho^*(\alpha)$.
\end{theorem}

We restricted our attention to the case of strong sparsity, that is,
$\alpha> 1/2$, as we may cover the whole range $\alpha\in(0,1]$ by
combining the ANOVA and the higher criticism tests (with a simple
Bonferroni correction), obtaining an adaptive test operating under
weaker constraints on the coherence $\gamma$.
That said, we mention that the higher criticism test is near-optimal in
the setting of Theorem \ref{thmSp-ub2} when, under the alternative, the
nonzero coefficients are not too spread out (restriction on the dynamic
range) and the amplitude is sufficiently large. This is the case, for
instance, when all nonzero coefficients are equal to $A$ in absolute
value with $A^2 S/\sqrt{p} > p^\eta$ for some $\eta> 0$ fixed.

The paper \cite{ingster2010detection} studies three tests assuming a
random design X. The first is based on $\|\by\|^2$ and is studied in
the nonsparse case where $S = p$, whereas the second is based on $\|\bX
^T \by\|^2$. The combined test is very similar to ANOVA and the
authors obtain the equivalent of Proposition \ref{prpmoderate-ub} for random
design matrices $\bX$ having standardized independent entries with
uniformly bounded fourth moment. Reference~\cite{ingster2010detection} also
considers the test based on the higher criticism applied to $|\vct{x}_j^T \by|/\|\by\|$ and the equivalent
of Theorems \ref{thmSp-ub} and \ref{thmSp-ub2} are established
under the assumption that the design matrix $\bX$ has i.i.d. standard
normal entries. Averaging over a
random design $\bX$ with standardized independent entries effectively
reduces to an orthogonal design, resulting in much weaker (implicit)
assumptions; no randomness assumptions on $\beta$---since this
randomness is carried by $\bX$---and no discretization of the
thresholds in the higher criticism statistic. In stark contrast, we
consider the design fixed (although it can of course be generated in a
random fashion).


Turning our attention to the Max test now, the results available for
orthogonal designs remain valid under similar conditions on the matrix
$\bX$.

\begin{theorem} \label{thmmax} Let $S = p^{1 - \alpha}$ and assume
that $\bX^T \bX\in\cS_p(\gamma, \Delta)$ with the
following parameter asymptotics: (1) $\Delta= O(p^{\eps})$, for all
$\eps> 0$ and (2) $\gamma^2 p^{1-\alpha}\times\break (\log p)^3 \to0$.
\begin{itemize}
\item In SFEM, the Max test is asymptotically powerful if $A \geq
\sqrt{2 r \log p}$ with $r > \rho_{\mathrm{Max}}(\alpha)$, and
asymptotically powerless if $r < \rho_{\mathrm{Max}}(\alpha)$.
\item In SREM, the Max test is asymptotically powerful for a fixed
signal level obeying $\tau> \rho^*_{\mathrm{rand}}(\alpha)$, and
asymptotically powerless if $\tau< \rho^*_{\mathrm{rand}}(\alpha)$.
\end{itemize}
%
The above holds for all $\alpha\in(1/2, 1]$.
\end{theorem}

This theorem justifies the assertion made in the \hyperref
[secintro]{Introduction}, which
stated that one could detect a linear relationship between the
response and a few covariates even though those covariates that were
mostly correlated with the response were not in the model. To
clarify, consider SFEM and $\alpha\in(1/2, 3/4]$. Then, for $A =
\sqrt{2 r \log p}$ with $\rho^*(\alpha) < r < \rho_{\mathrm{Max}}(\alpha)$,
the Max test is asymptotically powerless, whereas the test based on
$H^*$ has full power asymptotically. In particular, in the regime in
which the Max test is powerless, with high probability the entry of
$\bX^T y$ which achieves the maximal magnitude corresponds to a covariate
not in the support of $\beta$. (This is explicitly demonstrated in the
proof of Theorem \ref{thmmax}.) In the proof, we use fine asymptotic results
for the maximum of correlated normal random variables due to Berman
\cite{MR0161365} and Deo \cite{deo1972some}.

We pause here to comment on the situation in which the variance of the
noise (denoted $\sigma^2$) is unknown and must be estimated. As for
the identity design, the results in this section hold with $\by$
replaced by $\by/\hat{\sigma}$ with the proviso that $\hat{\sigma
}$ is
any accurate estimate with a slight upward bias to control the
significance level. Formally, suppose we have an estimator obeying
%
%
\begin{equation}
\label{eqconsistent}
\mathbb{P}\bigl(\sigma\le\hat\sigma\le(1+a_n) \sigma\bigr) \goto1
\end{equation}
and
{$a_n p^{1/2 - \epsilon} \goto0$ for all $\epsilon> 0$.}
We would then apply our methodology to $\by/\hat{\sigma}$. On the one
hand, it follows from the monotonicity of our statistic that the
asymptotic probability of type I errors is no worse than in the case of
known variance since we use an estimate which is biased upward. On the
other hand, consider an alternative with $S = p^{1-\alpha}$ and
amplitudes set to $A = \sigma\sqrt{2 r \log p}$, $r > \rho^*(\alpha)$.
The gap between $r$ and $\rho^*(\alpha)$ is sufficient to reject the
null. Indeed, $H^*$ is applied to~$\by/\hat{\sigma}$, leading to a
normalized amplitude equal to $\sqrt{2 r' \log p}$, where $r' :=
(\sigma/\hat{\sigma})^2 r$ is greater than $\rho^*(\alpha)$ in the
limit. (The contribution over the complement of the support of $\beta$
is negligible because $\hat{\sigma} -\sigma$ is sufficiently small,
{and this is why we require $a_n p^{1/2 - \epsilon} \goto0$}.) The same
arguments apply to the ANOVA $F$-test and the Max test. We mention that
Hall and Jin \cite{hj09} discuss the same issue for the case of an
orthogonal design and colored noise with a covariance that may be
unknown. Note that \cite{ingster2010detection} treats the case of
unknown variance in detail when the design matrix $\bX$ has i.i.d.
standard normal entries.

We now discuss strategies for constructing estimators obeying
(\ref{eqconsistent}). There are many possibilities and we choose to
discuss a simple estimate applying in the case of strong sparsity
$\alpha\in(1/2, 1]$, where signals are near the detection boundary,
so that $\|X \beta\|^2 /(\sigma^2 \sqrt{n}) \rightarrow0$ (this is
the interesting regime).
{For concreteness, assume that $n < p = O(n^{1 + \epsilon})$ for
all $\epsilon> 0$.}
As noted in Section
\ref{secpriorWork}, $\| \vct{y}\|^2/\sigma^2$ has the chi-square
distribution with $n$ degrees of freedom and noncentrality parameter
$\|\bX\bbeta\|^2/\sigma^2$, and, thus,
\[
\mathbb{P}\bigl(\sigma\bigl(1-s_n/\sqrt{n}\bigr) \le\|\by\|/\sqrt{n} \le\sigma\bigl(1 +
s_n/\sqrt{n}\bigr)\bigr) \goto1
\]
as long as $s_n \rightarrow\infty$. Now let $t_n \to\infty$ slowly
(say, $t_n = \log n$) and define $\hat{\sigma} := \|\by\|(1/\sqrt
{n} +
t_n/n)$. This estimator obeys (\ref{eqconsistent}).

\subsection{Normal designs}

A common assumption in multivariate statistics is that the rows of the
design matrix are independent draws from the multivariate normal
distribution $\mathcal{N}(0,{\bolds\Sigma})$. Our results apply provided
that ${\bolds\Sigma}$ obeys the assumptions about $\bX^T\bX$.
\begin{cor} \label{corstat} Suppose the rows of $\bX$ are independent
samples from $\mathcal{N}(0,{\bolds\Sigma})$, and ${\bolds\Sigma} \in
\cS_p(\gamma,\Delta)$ (the columns are normalized). Then
the conclusions of Theorems \ref{thmSp-lb}, \ref{thmSp-ub} and
\ref{thmmax} are all valid, provided that $\sqrt{n^{-1} \log
p}$ obeys the conditions imposed on $\gamma$.
\end{cor}

We remark that if the columns are not normalized so that the rows of
$\bX$ are independent samples from $\mathcal{N}(0,{\bolds\Sigma})$, the
same result holds with a threshold $A$ replaced by $A/\sqrt{n}$. This
holds because the norm of each column is sharply concentrated around
$\sqrt{n}$.

\section{Some special designs}
\label{secfull}

We consider correlation matrices which have a substantial portion
of large entries. In general, the detection threshold may depend upon
some fine details of $\bX$, but we give here some representative results
applying to situations of interest.

We first examine the simple, yet important and useful example of
constant correlation, where $\vct{x}_j^T \vct{x}_k = 1$ if $j = k$, and $=
\gamma$ if $j \neq k$.\footnote{Whether such a family of vectors
exists for special values of $\gamma$ is a nontrivial matter, and we
refer the reader to the literature on equiangular lines; see
\cite{MR0307969}, for example.} We impose $0 < \gamma< 1$ to make
sure that $\bX^T \bX$ is at least positive definite as $p \to\infty$
(this implies that $\bX^T \bX$ has full rank which in turn imposes $p
\leq n$). The balanced one-way design has this structure since it can
be modeled by the matrix
\[
\bX= \frac{1}{\sqrt{2k}} \left[\matrix{
{\mathbf1} & {\mathbf0} & \cdots& {\mathbf0} \cr
{\mathbf0} & {\mathbf1} & \cdots& {\mathbf0} \cr
\vdots& \vdots& \vdots& \vdots\cr
{\mathbf0} & {\mathbf0} & \cdots& {\mathbf1} \cr
- {\mathbf1} & - {\mathbf1} & \cdots& - {\mathbf1}}\right],
\]
where each vector in this block representation is $k$-dimensional.
Without further assumptions on $\bbeta$, this design is equivalent to
(\ref{one-way}) with the constraint ${\mathbf1}^T \bbeta= 0$, except for
the normalization. With this definition, $\bX^T\bX$ has diagonal
entries equal to $1$ and off-diagonal entries equal to $1/2$ so we are
in the setting---with $\gamma= 1/2$---of our next result below.
\begin{theorem} \label{thmzeta} Suppose that $\vct{x}_j^T \vct{x}_k$
is equal to
$1$ if $j = k$ and $\gamma$ otherwise, and that the sparsity exponent
obeys $\alpha\in(1/2,1]$. Then without further assumption, the
conclusions of Theorems \ref{thmSp-lb}, \ref{thmSp-ub} and
\ref{thmmax} remain valid with the bounds on $A$ and $\tau$
divided by $\sqrt{1-\gamma}$.
\end{theorem}

The balanced, one-way design may be seen either as an orthogonal
design with a linear constraint, or a constant-correlation design
without any constraint. More generally, a multi-way design is easily
defined as an orthogonal design with a set of linear constraints.
Specifically, suppose the coordinates of $\bbeta$ are indexed by an
$m$-dimensional index vector, so that
\[
\bbeta= \bigl(\beta_\vct{j}\dvtx\vct{j} = (j_1,\ldots, j_m), j_s \in
[p_s]\bigr),\qquad p = \prod_{s=1}^m p_s.
\]
We assume the design is balanced with $k$ replicates per cell so that
$n = p k$. With any fixed order on the index set, say, the
lexicographic order, the design matrix is the same as in the balanced,
one-way design (\ref{one-way}). Here, $\bbeta$ obeys the linear
constraints
%
%
\begin{equation}
\label{constraints}
\sum_{s \neq t} \sum_{j_s = 1}^{p_s} \beta_{j_1 \cdots j_m} = 0
\end{equation}
for all $j_t \in[p_t]$ and $t \in[m]$ (there are $\sum_{t=1}^m p_t$
constraints). As in the balanced, one-way design, Theorem \ref{thmortho}
applies to the balanced, multi-way design. The
argument for the lower bound is at the end of the proof of
Theorem \ref{thmSp-lb}. The proof of the upper bounds is exactly as in the
case of any other orthogonal design. Finally, embedding the linear
constraints into the design matrix leads to a family of designs with a
``full'' correlation structure with off-diagonal elements which, in
general, are not of the same magnitude unless the design is one-way.

\section{Numerical experiments}
\label{secnumerics}

We complement our study with some numerical simulations
which illustrate the empirical performance for finite sample
sizes. Here, $\bX$ is an $n \times p$ Gaussian design with
i.i.d. standard normal entries, and normalized columns. We study fixed
effects and investigate the performance of ANOVA, the higher criticism\footnote{We do not use the discretization here.} and the Max
test. We also compare the detection limits with those available in the
case of the $p \times p$ identity design, since the theory developed
in Corollary \ref{corstat} predicts that the detection boundaries are
asymptotically identical (provided $n$ grows sufficiently rapidly).

%
%
\begin{figure}

\includegraphics{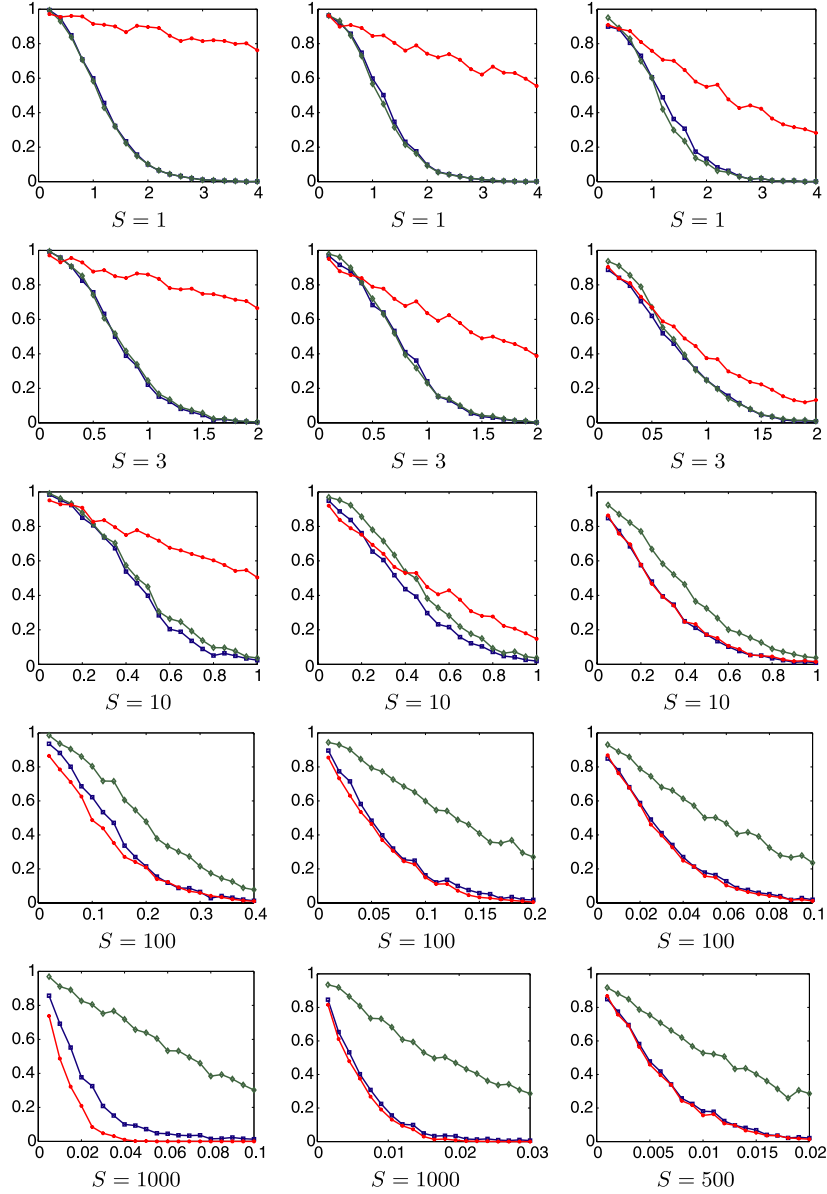}

\caption{Left column: identity design with $p = 10\mbox{,} 000$. Middle
column: Gaussian design with $p = 10\mbox{,} 000$ and $n = 2\mbox{,} 000$. Right
column: Gaussian design with $p = 10\mbox{,} 000$ and $n = 500$. Sparsity
level $S$ is indicated below each plot. In each plot, the empirical
risk (based on $1\mbox{,} 000$ trials) of each method [ANOVA (red bullets);
higher criticism (blue squares); Max test (green diamonds)] is
plotted against $r$ (note the different scales).}
\label{figmoderate}
\end{figure}

\begin{figure}

\includegraphics{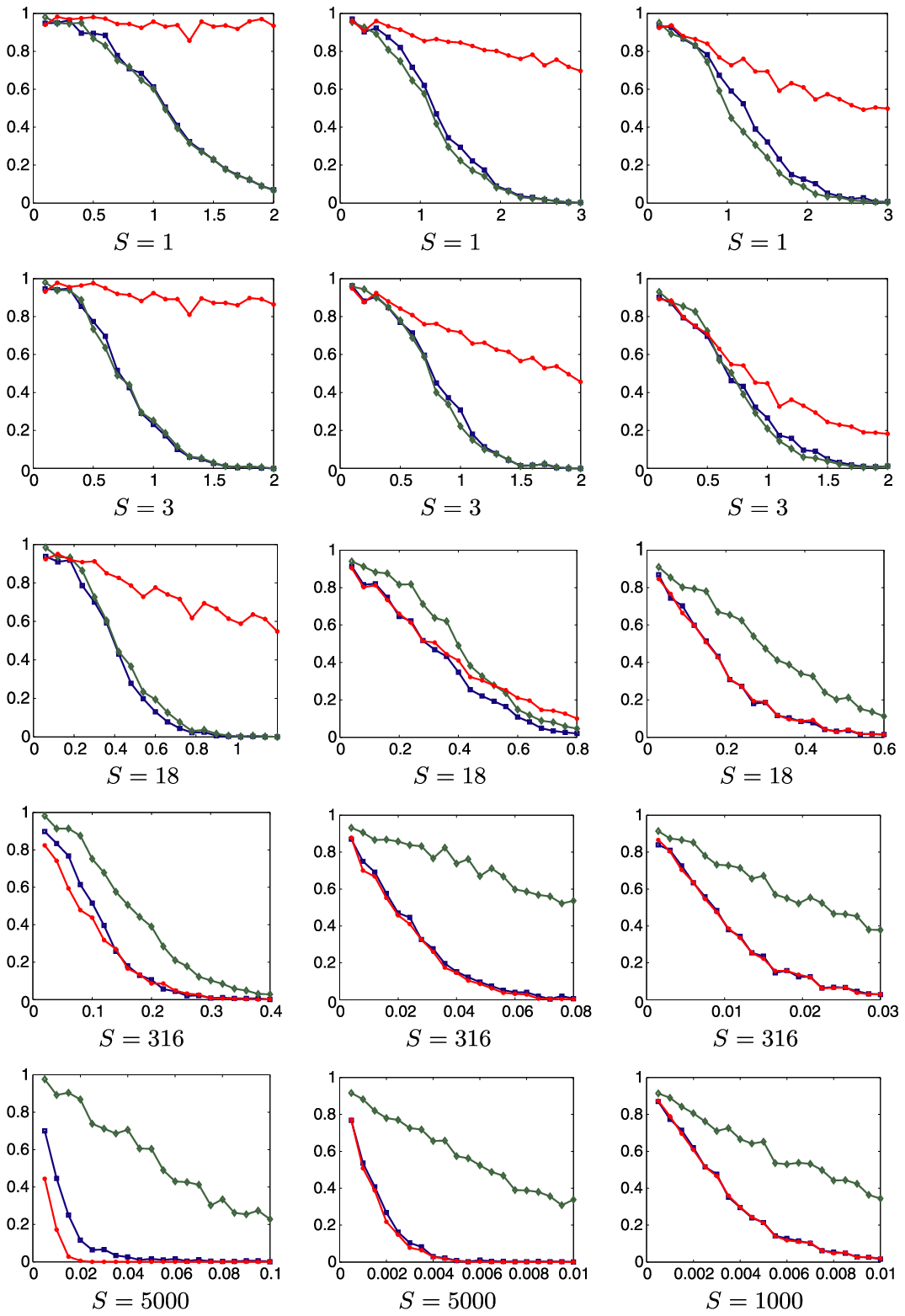}

\caption{Left column: identity design with $p = 100\mbox{,} 000$. Middle
column: Gaussian design with $p = 100\mbox{,} 000$ and $n = 5\mbox{,} 000$. Right
column: Gaussian design with $p = 100\mbox{,} 000$ and $n =
1\mbox{,} 000$. Sparsity level $S$ is indicated below each plot. In each
plot, the empirical risk (based on 500 trials) of each method [ANOVA
(red bullets); higher criticism (blue squares); Max test (green
diamonds)] is plotted against $r$ (note the different scales).}
\label{figlarge}
\end{figure}

We performed simulations with matrices of sizes $500 \times10\mbox{,} 000$,
$2\mbox{,} 000 \times10\mbox{,} 000$, $1\mbox{,} 000 \times100\mbox{,} 000$ and $5\mbox{,} 000
\times100\mbox{,} 000$, various sparsity levels, and strategically selected
values of $r$. Each data point corresponds to an average over
$1\mbox{,} 000$ trials in the case where $p = 10\mbox{,} 000$, and over $500$
trials when $p = 100\mbox{,} 000$. A new design matrix is sampled for each
trial. The performance of each of the three methods is computed in
terms of its best (empirical) risk defined as the sum of probabilities
of type I and II errors achievable across all thresholds. The results
are reported in Figures \ref{figmoderate} and \ref{figlarge}. As
expected, the detection thresholds for the Gaussian design are quite
%
%
close to those available for the identity design. The performance of
ANOVA improves very quickly as the sparsity decreases, dominating the
Max test with $S = \sqrt{p}$; its performance also improves as $n$
becomes smaller, in accordance with~(\ref{anova}). The performance of
the Max test follows the opposite pattern, degrading as $S$ increases.
Interestingly, the higher criticism remains competitive across the
different sparsity levels.

\section{Discussion}
\label{secdiscussion}

It is possible to extend our results to setups with correlated errors,
with known covariance. As discussed in Section \ref{secintro}, suppose
$\bz$ in (\ref{setting}) is $\cN({\mathbf0}, \mtx{V})$. We may then
whiten the noise by multiplying both sides of (\ref{setting}) by $\mtx
{L}^{-1}$, where $\mtx{L} \mtx{L}^T$ is a Cholesky decomposition of
$\mtx{V}$. This leads to a model of the form
\[
\by= \mtx{L}^{-1} \bX\bbeta+ \bz,
\]
which is our problem with $\mtx{L}^{-1} \bX$ instead of $\bX$. In
some situations, the noise covariance matrix may not be known and we
refer to \cite{hj09} for a brief discussion of this issue.

Although several generalizations are possible, an interesting open
problem is to determine the detection boundary for a given sequence of
designs $\{\bX_{n \times p}\}$ with $n$ and $p$ growing to
infinity. We have seen that if most of the predictor variables are
only weakly correlated, then the detection boundary is as if the
predictors were orthogonal. Similar conclusions for certain types of
square designs in which $n = p$ are also presented in the work of Hall
and Jin \cite{hj09}. Although we introduced some sharp results in
Section \ref{secfull} corresponding to some important design matrices, the
class of matrices for which we have definitive answers is still quite
limited. We hope other researchers will engage this area of research
and develop results toward a general theory.



\section*{Acknowledgments}

We would like to thank Chiara Sabatti for stimulating discussions and
for suggesting improvements on an earlier version of the manuscript,
and Ewout van den Berg for help with the simulations. We also thank
the anonymous referees for their inspiring comments which helped
us improve the content of the paper.

\begin{supplement}[id=suppA]
\stitle{Supplement to ``Global testing under sparse alternatives: ANOVA,
multiple comparisons and the higher criticism''}
\slink[doi]{10.1214/11-AOS910SUPP} 
\sdatatype{.pdf}
\sdescription{In the supplement, we prove the results stated in the
paper. Though the method of proof has the same structure as the
corresponding situation in the classical setting with identity design
matrix, extra care is required to deal with dependencies.}
\end{supplement}

%

\printaddresses

\end{document}